\documentclass[11pt]{article}
\usepackage{amssymb,amsfonts,amsmath,latexsym,epsf,tikz,url}

\newtheorem{theorem}{Theorem}[section]

\newtheorem{corollary}[theorem]{Corollary}

\newcommand{\proof}{\noindent{\bf Proof.\ }}
\newcommand{\qed}{\hfill $\square$\medskip}

\begin{document}

\title{The distinguishing index of graphs with at least one cycle is not more than its distinguishing number} 

\author{
Saeid Alikhani\footnote{Corresponding author} \and Samaneh Soltani 
}

\date{\today}

\maketitle

\begin{center}
Department of Mathematics, Yazd University, 89195-741, Yazd, Iran\\
{\tt alikhani@yazd.ac.ir, s.soltani1979@gmail.com }

\end{center}

\begin{abstract}
The distinguishing number (index) $D(G)$ ($D'(G)$) of a graph $G$ is the least integer $d$ such that $G$ has an vertex (edge) labeling with $d$ labels that is preserved only by the trivial automorphism. It is known that for every graph $G$ we have $D'(G) \leq D(G) + 1$.  The complete characterization of finite trees $T$ with $D'(T)=D(T)+ 1$ has been given recently.  
 In this note we show that if $G$ is a finite connected graph with at least one cycle, then $D'(G)\leq D(G)$. Finally, we characterize all connected graphs for which $D'(G) \leq D(G)$.  
\end{abstract}


\noindent
{\bf AMS Subj.\ Class.\ (2010)}: 05C15, 05E18

\section{Introduction}
\label{sec:intro}

Let $G = (V(G), E(G))$ be a finite graph and  ${\rm Aut}(G)$ be its automorphism group. 
 A labeling $\phi: V(G) \rightarrow [r]$ is {\em distinguishing} if no non-trivial element of ${\rm Aut}(G)$ preserves all the labels; such a labeling $\phi$ is a {\em distinguishing $r$-labeling}. More formally, $\phi$ is a distinguishing labeling if for every $\alpha\in {\rm Aut}(G)$, $\alpha\ne {\rm id}$, there exists $x\in V(G)$ such that $\phi(x) \neq \phi(\alpha(x))$. The {\em distinguishing number} $D(G)$ of a graph $G$ is the smallest $r$ such that $G$ admits a distinguishing $r$-labeling. 

The introduction of the distinguishing number in 1996 by Albertson and Collins~\cite{albertson-1996} was a great success;   by  now  more than   one  hundred  papers  were  written  motivated  by  this  paper!  The distinguishing number of several families of graphs, including trees, hypercubes, generalized Petersen graphs, friendship and book graphs have been studied in \cite{alikhani-2017+,Bogstad,Chan,cheng-2006}, and \cite{Potanka}. 

Similar to this definition, Kalinowski and Pil\'sniak~\cite{kalinowski-2015} have defined the distinguishing index $D'(G)$ of $G$. The {\em distinguishing index} $D'(G)$ of a graph $G$ is the smallest integer $d$ such that $G$ has an edge labeling with $d$ labels that is preserved only by the trivial automorphism. 
Generally $D'(G)$ can be arbitrary smaller than $D(G)$, for instance if $p\geq 6$, then $D'(K_p)=2$ and $D(K_p)=p$. Conversely, there is an upper bound on $D'(G)$ in terms of $D(G)$.
In~\cite[Theorem 11]{kalinowski-2015} (see also~\cite[Theorem 8]{lehner-2017+} for an alternative proof) it is proved that if $G$ is a connected graph of order at least $3$, then
\begin{equation*}
D'(G) \leq D(G) + 1.
\end{equation*}

Authors in \cite{Sandi etal-2017} gave  complete characterization of finite trees with $D' = D + 1$. They  defined a family ${\cal T}$ as follows. It consists of those trees $T$ of order at least $3$, for which the following conditions are fulfilled. 
\begin{enumerate}
\item The tree  $T$ is a bicentric tree with the central edge $e=vw$. 
\item There is an isomorphism between the rooted trees $T_v$ and $T_w$ (Note that $T_v$ and $T_w$ are components of $T-e$, where $e=vw$, $v\in T_v$ and $w \in T_w$). 
\item There is a unique distinguishing edge-labeling of the rooted tree $T_v$ using $D(T)$ labels.
\end{enumerate}
The following theorem  states that the family ${\cal T}$ contains all trees with $D'(T) = D(T)+1$.
\begin{theorem}{\rm \cite{Sandi etal-2017}}
\label{thm:trees}
Let $T$ be a tree of order at least $3$. Then 
$$D'(T) = \begin{cases}
  D(T) + 1 & {\rm if }\;T\in {\cal T}, \\
  D(T) & {\rm otherwise}.
\end{cases}$$
\end{theorem}

They also  showed that if $G$ is a connected unicyclic graph, then $D'(G) = D(G)$, showing that the inequality is never sharp for unicyclic graphs.
We shall  show that the distinguishing index of connected graphs  with at least one cycle is at most equal with its  distinguishing number.
We need to recall about ``orbit of a set". Let $\Gamma$ be a group acting on a set $X$. If $g$ is an element of $\Gamma$ and $x$ is in $X$, then denote the action of $g$ on $x$ by $g(x)$. The action of $\Gamma$ on $X$  is {\em faithful} if for any two distinct elements $g, h \in \Gamma$ there exists an element $x\in X$ such that $g(x)\neq h(x)$. Write ${\rm Orbit}(x)$ for the orbit containing $x$, i.e., ${\rm Orbit}(x)=\{g(x):~ g\in \Gamma\}$. Let $Y=\{y_1,\ldots, y_k\}$ be a subset of $X$. The mean of  orbit of the subset $Y$ of $X$ is as follows. 
\begin{equation*}
{\rm Orbit}(Y)=\big\{\{g(y_1),\ldots, g(y_k)\}:~ g\in \Gamma \big\}.
\end{equation*}

If $G$ is  a graph, then  ${\rm Aut}(G)$ acts faithfully on $V(G)$.  As an example of orbit of a set, consider the complete graph $K_n$ with vertex set $\{v_1,\ldots , v_n\}$. It is clear that there exists exactly one cycle  of length $n$ with alternative vertices $v_1,\ldots , v_n$ in $K_n$, say $C$, so ${\rm Orbit}(C)=\{C\}$. If $C'$ is an $i$-cycle in $K_n$, where $i< n$, then $|{\rm Orbit}(C')|\geq 2$, since there are at least two  $i$-cycle in $K_n$ for any  $i< n$.

\medskip
In this paper we show that the distinguishing index of finite, connected graphs  with at least one cycle is at most equal with the distinguishing number. Finally, we give a characterization of all finite,  connected graphs  with $D'=D+1$, by Theorem \ref{thm:trees}.

\section{Main result}



\begin{figure}[ht]
\begin{center}
\includegraphics[width=0.8\textwidth]{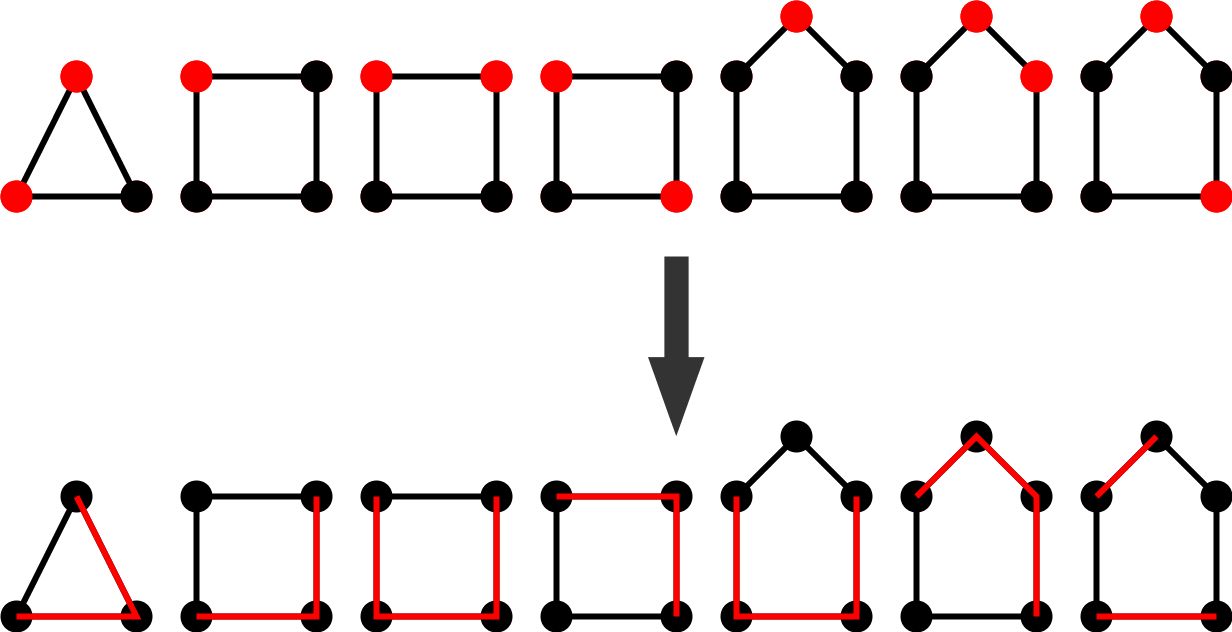}
\caption{All non-equivalent vertex $2$-labelings of $C$ and their respective transformations to edge $2$-labelings of $C$}
\label{fig:vertex-to-edge}  	
\end{center}
\end{figure}

\begin{theorem}\label{Maintheorem}
If  $G$ is a finite, connected graph with at least one cycle, then $D'(G)\leq D(G)$.
\end{theorem}
\proof
 Let $C_0$ be a cycle in $G$ which has the smallest orbit among other cycles in $G$ (note that the orbit of the cycle $C_0$ is a set contains all cycles in $G$ such that the cycle $C_0$ can be mapped on those under the automorphisms of $G$). Without loss of generality we  assume that  ${\rm Orbit}(C_0)=\{C_0,C_1, \ldots , C_k\}$, and  hence $|{\rm Orbit}(C)|\geq k$ for any cycle $C$ in $G$.
We want to present an edge distinguishing labeling of $G$, say $L$,  with at most $D(G)$ labels. For this purpose, we consider a vertex distinguishing labeling of $G$ with $D(G)$ labels $\{1,\ldots , D(G)\}$, say $\varphi$, and continue by the following steps:
\begin{enumerate}
\item[Step 1)] Similar to the method used in Figure \ref{fig:vertex-to-edge}, we can consider an edge labeling of $C_i$, $0\leq i \leq k$, using the labels of vertices $C_i$ in $\varphi$, such that
\begin{enumerate}
\item If $f$ is an automorphism of $G$ mapping the cycle $C_i$ to itself and preserving the edge labeling of $C_i$, then $f$ preserves the restricted vertex labeling $\varphi$ to $C_i$.

In fact, if $f$ is an automorphism of $G$ mapping the cycle $C_i$ to itself, then $f$ is an automorphism of cycle graph $C_i$. Hence, by a method similar to Figure \ref{fig:vertex-to-edge}, we can label the edges of $C_i$, using the labels of vertices $C_i$ in $\varphi$, such that if $f$ preserves the edge labeling of $C_i$, then $f$ preserves the restricted vertex labeling $\varphi$ to $C_i$.
\item If $f$ is an automorphism of $G$ mapping the cycle $C_i$ to $C_j$ where $i\neq j$ and $i,j\in \{0,1,\ldots , k\}$, and preserving the edge labeling between them, then $f$ preserves the restricted vertex labeling $\varphi$ to $C_i \cup C_j$.

In fact, if $f$ is an automorphism of $G$ mapping the cycle $C_i$ to $C_j$ where $i\neq j$ and $i,j\in \{0,1,\ldots , k\}$, then $f$ is an automorphism of $C_i$ to $C_j$. Hence we can label the edges $C_i$ and $C_j$, using the labels of vertices $C_i$ and $C_j$ in $\varphi$ respectively,  such that if $f$ preserves the edge labeling between them, then $f$ preserves the restricted vertex labeling $\varphi$ to $C_i \cup C_j$.
\end{enumerate}
\item[Step 2)] For labeling the remaining unlabeled edges of $G$, we denote the subgraph $\bigcup_{i=0}^k C_i$ by $H$, and use the following formula for any $xy\in E(G)\setminus E(H)$:
\begin{equation*}
L(xy)=\left\{
\begin{array}{ll}
\varphi (x) & \text{if}~~{\rm dist}_G (x,H) > {\rm dist}_G (y,H),\\
\varphi (y) & \text{if}~~{\rm dist}_G (x,H) < {\rm dist}_G (y,H),\\
1 & \text{if}~~{\rm dist}_G (x,H) = {\rm dist}_G (y,H).
\end{array}\right.
\end{equation*} 
\end{enumerate}
We claim that the edge labeling $L$ is distinguishing. In fact, if $f$ is an automorphism of $G$ preserving this edge labeling, then since $f$ maps the subgraph $H$ to itself and since $f$ preserves the edge labeling of $H$, so $f$ preserves the vertex labeling of $H$, by Step 1. On the other hand, since $f$ preserves the  distances from vertices of $G$ to $H$, so we can conclude that $f$ preserves the vertex labeling of $G$, i.e., $\varphi$, and hence $f$ is the identity automorphism.\qed

\begin{figure}[ht]
\begin{center}
\includegraphics[width=0.8\textwidth]{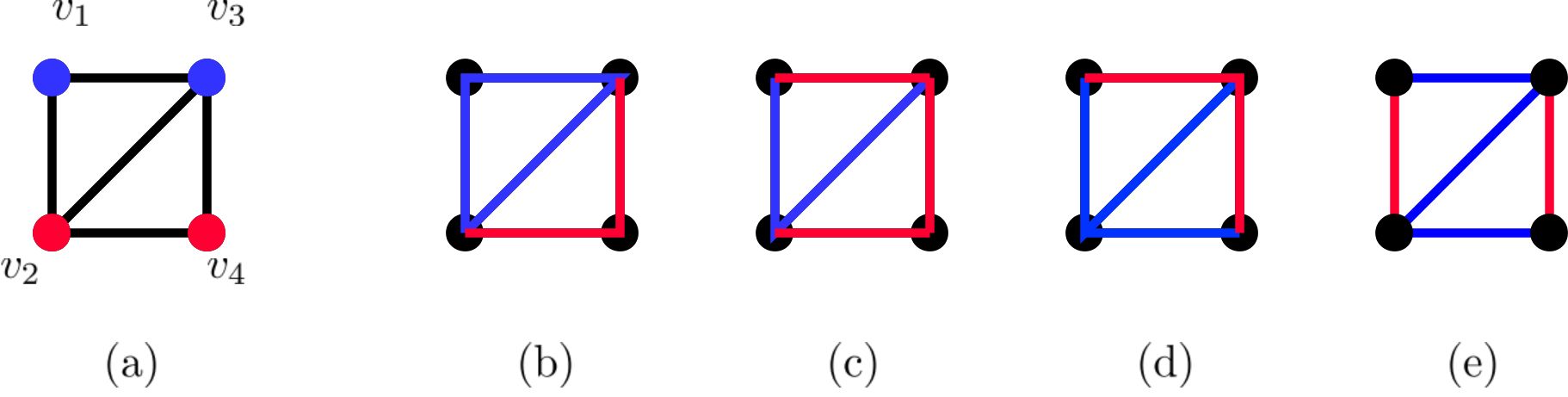}
\caption{\label{55}  Colored $K_4-e$.}
\end{center}
\end{figure}

 \noindent{\bf Remark.} In the  proof of Theorem \ref{Maintheorem}, it seems that the selection of the smallest orbit cycle is not used at all, but has been used indirectly. It is possible that there exist small cycles inside large cycles. For example, consider the graph $K_4 -e$. In this case, we have two 3-cycles and exactly one 4-cycle inside $K_4-e$ such that 3-cycles have been lie in the 4-cycle. By Figure \ref{55} we  show that Step 1 of proof of Theorem \ref{Maintheorem} is impossible for 3-cycles.
Figure (a) shows a 2-distinguishing labeling of graph.  Since edge $\{v_2,v_3\}$ is fixed under each automorphism of graph, so the label of that can be arbitrary. Figures (b), (c), (d) and (e) are all possible non-isomorphic colorings of graph, without respect to the label of edge $\{v_2,v_3\}$. Thus Step 1 (a) of proof of Theorem \ref{Maintheorem} is not satisfied.
 \begin{itemize}
\item Figure (b): There exist automorphism $f$ with $f(v_1)=v_1$, $f(v_4)=v_4$, $f(v_2)=v_3$ and $f(v_3)=v_2$ such that $f$ preserves the edge labeling of each of 3-cycles, but $f$ does not preserve the vertex labeling of them. Thus Step 1 (b) of proof of Theorem \ref{Maintheorem} is not satisfied.
\item  Figure (c):  There exist automorphism $f$ with $f(v_1)=v_1$, $f(v_4)=v_4$, $f(v_2)=v_3$ and $f(v_3)=v_2$ such that $f$ preserves the edge labeling of 3-cycle $v_2v_3v_4$, but $f$ does not preserve the vertex labeling of this 3-cylces. Thus Step 1 (a) of proof of Theorem \ref{Maintheorem} is not satisfied.
\item  Figure (d): There exist automorphism $f$ with $f(v_1)=v_4$, $f(v_4)=v_1$, $f(v_2)=v_2$ and $f(v_3)=v_3$ such that $f$ preserves the edge labeling of between 3-cycles $v_4v_2v_3$ and $v_1v_2v_3$, but $f$ does not preserve the vertex labeling between them. Thus Step 1 (b) of proof of Theorem \ref{Maintheorem} is not satisfied.
\item  Figure (e): There exist automorphism $f$ with $f(v_1)=v_4$, $f(v_4)=v_1$, $f(v_2)=v_3$ and $f(v_3)=v_2$ such that $f$ preserves the edge labeling of between 3-cycles $v_4v_2v_3$ and $v_1v_2v_3$, but $f$ does not preserve the vertex labeling between them. Thus Step 1 (b) of proof of Theorem \ref{Maintheorem} is not satisfied.
 \end{itemize}

 We end the paper with  a characterization of  finite, connected graphs  with $D'=D+1$ which follows from  Theorems \ref{thm:trees} and \ref{Maintheorem}. 
\begin{corollary}
Let $G$ be a finite, connected graph of order $n\geq 3$. Then $D'(G)=D(G)+1$ if and only if $G$ is a tree in ${\cal T}$.
\end{corollary}

\end{document}